\documentclass[twoside,12pt,leqno]{article}
\usepackage{amsmath, amsthm, amsfonts, amssymb, color}
\usepackage{mathrsfs}
\usepackage{fancyhdr}

\newcommand{\be}{\begin{eqnarray}}
\newcommand{\ee}{\end{eqnarray}}
\newcommand{\ce}{\begin{eqnarray*}}
\newcommand{\de}{\end{eqnarray*}}

\newtheorem{thm}{Theorem}[section]

\theoremstyle{definition}

\definecolor{wco}{rgb}{0.5,0.2,0.3}

\numberwithin{equation}{section}
\theoremstyle{remark}
\newtheorem{rem}{Remark}[section]

\newcommand{\ua}{\uparrow}

\title{{\bf Harnack Inequality and Applications for Stochastic
 Evolution Equations with Monotone Drifts}
\footnote{Supported in part by the DFG through the Internationales
 Graduiertenkolleg
``Stochastics and Real World Models'', the SFB 701 and NNSFC(10721091).} }
\author{{\bf Wei Liu }\\
\footnotesize{School of Mathematical Sciences, Beijing Normal
University, 100875 Beijing, China}\\
{\footnotesize Fakult\"at f\"ur Mathematik, Universit\"at Bielefeld,
D-33501 Bielefeld, Germany
 }\\
\footnotesize{E-mail: wei.liu@uni-bielefeld.de}}

\date{}
\begin{document}
\maketitle

\def\R{\mathbf R}  \def\ff{\frac} \def\ss{\sqrt} \def\BB{\mathbf
B}
\def\N{\mathbf N} \def\kk{\kappa} \def\m{{\mu}}
\def\dd{\delta} \def\DD{\Delta} \def\vv{\varepsilon} \def\rr{\rho}
\def\<{\langle} \def\>{\rangle} \def\GG{\Gamma} \def\gg{\gamma}
  \def\nn{\nabla} \def\pp{\partial} \def\tt{\tilde}
\def\d{\text{\rm{d}}} \def\bb{\beta} \def\aa{\alpha} \def\D{\mathcal D}
\def\E{\mathbf E} \def\si{\sigma} \def\ess{\text{\rm{ess}}}
\def\beg{\begin} \def\beq{\begin{equation}}  \def\F{\mathcal F}
\def\Ric{\text{\rm{Ric}}} \def\Hess{\text{\rm{Hess}}}\def\B{\mathbf B}
\def\e{\text{\rm{e}}} \def\ua{\underline a} \def\OO{\Omega}
 \def\b{\mathbf b}
\def\oo{\omega}     \def\tt{\tilde} \def\Ric{\text{\rm{Ric}}}
\def\cut{\text{\rm{cut}}} \def\P{\mathbf P} \def\ifn{I_n(f^{\bigotimes
 n})}
\def\fff{f(x_1)\dots f(x_n)} \def\ifm{I_m(g^{\bigotimes m})}
 \def\ee{\varepsilon}
\def\pm{\pi_{{\bf m}}}   \def\p{\mathbf{p}}   \def\ml{\mathbf{L}}
 \def\C{\mathcal C}      \def\aaa{\mathbf{r}}     \def\r{r}
\def\gap{\text{\rm{gap}}} \def\prr{\pi_{{\bf m},\varrho}}
  \def\r{\mathbf r}
\def\Z{\mathbf Z} \def\vrr{\varrho} \def\ll{\lambda}
\def\L{\mathcal L}

\def\[{{\Big[}}
\def\]{{\Big]}}

\def\({{\Big(}}
\def\){{\Big)}}

\def\xy{{\|X_t-Y_t\|}}
\def\x{{\|X_t\|}}
\def\y{{\|Y_t\|}}


\begin{abstract} As a  Generalization to \cite{W07} where the  dimension-free Harnack
inequality was established  for
stochastic  porous media equations, this paper presents analogous
results for a large class of stochastic evolution equations with
general monotone drifts. Some ergodicity, compactness
 and contractivity
properties are  established for the associated transition
semigroups. Moreover, the exponential convergence of the transition semigroups to
invariant measure and the existence of a spectral gap are also derived. As examples, the main results are applied to
many  concrete SPDEs such as stochastic reaction-diffusion
equations,  stochastic porous media equations and
 the stochastic $p$-Laplace
equation in Hilbert space.
\end{abstract}



\bigbreak

\section{Introduction and main results}

The dimension-free Harnack inequality has been a very efficient tool
for the study of diffusion semigroups in recent years. It was
 first introduced by Wang in
\cite{W97} for  diffusions on Riemannian manifolds, then this
infinite
 dimensional version of
Harnack inequality has been applied and extended intensively later, see
e.g. \cite{W99, W01, RW1, RW2} for applications to
 functional inequalities; \cite{AK,AZ,K} for
the study of short time behavior of infinite-dimensional
diffusions; \cite{GW04,W00} for the estimate of high order
eigenvalues,  and \cite{BGL} for applications to the
transportation-cost inequality and \cite{GW}
for heat kernel estimates.

Recently, the dimension-free Harnack inequality was established in
\cite{W07} for  stochastic  porous media
equations and in \cite{LW} for stochastic
 fast-diffusion
equations. As applications, the strong Feller property, estimates of
the transition density and some  contractivity properties were
obtained for the associated transition semigroups. The approach used
in \cite{LW,W07} is based on a coupling argument developed in
\cite{ATW}, where the Harnack inequality was studied for diffusion
semigroups on Riemannian manifolds with unbounded curvatures from below.
The advantage of this approach is that one can avoid the assumption that the
 curvature is lower bounded, which was used in previous works (cf.\cite{AK,AZ,BGL,RW1,RW2})
 in an essential way and would be very hard to verify in the
present framework of non-linear SPDE.

The aim of this paper is to establish
the analogous results for general stochastic evolution equations
within the variational  framework. More precisely, we mainly deal with
 stochastic evolution equations
with strongly dissipative drifts in Hilbert space,
whick cover many important
 types of SPDE such as  stochastic reaction-diffusion equations,
stochastic porous media
equations and the stochastic p-Laplace equation (cf.\cite{PR,KR,Zh}).
We first establish the Harnack inequality and the strong Feller property for the
associated transition semigroups, then it has been used to
derive some ergodicity and contractivity properties for the
 transition semigroups.
 In particular, we give a very easy proof for the (topological)
 irreducibility in Theorem \ref{T1.2} by using the
Harnack inequality.
Hence  the uniqueness of invariant
 measures  for the transition
semigroups is obtained without assuming  strict monotonicity for the drift,
which was required  in many
earlier works \cite{PR,W07,LW,RRW,DRRW}.  And we also derive the
 convergence rate of the transition semigroups to the invariant measure.
This implies a decay estimate of the solutions to the corresponding deterministic evolution equations (e.g.$p$-Laplace
equation, porous medium equation), which coincides with some well-known results
 in PDE theory. Moreover, some uniformly exponential ergodicity and
the existence of a spectral gap are also investigated.

Now we  describe our framework for SPDE in details. There
exist three main different approaches to analyze stochastic
partial
 differential equations in the literature. The martingale measure
 approach  was
 initiated by J.
 Walsh in \cite{Wa}. The  variational approach was first
 used by
 Pardoux \cite{Pa} to study SPDE, then this approach
 was further developed by Krylov and
 Rozovoskii  \cite{KR} and applied to non-linear filtering. Concerning
 the
 semigroup
 approach  we  refer to the classical book by Da Prato and Zabcyzk
 \cite{DaZa}. In this paper we will use the variational approach
 because we
  mainly treat nonlinear SPDE of evolutionary type.
 All kinds of dynamics with stochastic influence
in nature or man-made complex systems can be
 modeled by such equations. This type of SPDE
  has been
 studied intensively  in recent years, we refer to
 \cite{DRRW,GM,L08,RRW,Kim,PR,Zh} (and references therein)
 for various generalizations and applications.

Let $H$ be a separable Hilbert space with inner product
 $\<\cdot,\cdot\>_H$ and
$H^*$ its dual. Let  $V$ be a reflexive and separable Banach space
such
 that
 $V\subset H$  continuously and
densely. Then for its dual space $V^*$ it follows that $H^*\subset
V^*$
 continuously
and densely. Identifying $H$ and $H^*$ via the Riesz isomorphism we
 know that
$$V\subset H\equiv H^*\subset V^*$$
is a Gelfand triple.  If the dualization between $V^*$ and $V$ is
denoted by $_{V^*}\<\cdot,\cdot\>_V$ we have
$$ _{V^*}\<u,v\>_V= \<u,v\>_H \ \  \text{for all} \
u\in H, v\in V. $$

Suppose $W_t$ is a cylindrical Wiener process on a separable Hilbert
space $U$ w.r.t a complete filtered probability space
$(\Omega,\mathcal{F},\mathcal{F}_t,\mathbf{P})$,  and
$(L_2(U;H), \|\cdot\|_2)$ is the space of
all  Hilbert-Schmidt operators from $U$ to $H$. Now we consider the
following stochastic evolution  equation
 \beq
\label{1.1} \d X_t=A(t,X_t)\d t+B_t\d W_t,\ X_0=x\in H,
 \end{equation}
 where
$$A:[0,T]\times V\times\Omega\rightarrow V^*; \
B:[0,T]\times\Omega\rightarrow L_2(U,H)$$
are progressively measurable.
We first recall the classical result in \cite{KR}
for the existence and
uniqueness of strong solution.
For more general results we refer to \cite{GM,RRW,Zh}.

\beg{lem} {\bf (\cite{KR} Theorems II.2.1, II.2.2 )}
\label{L0}
 Consider the general stochastic evolution
 equation
 \begin{equation}\label{general equation}
\d X_t=A(t,X_t)\d t+B(t,X_t)\d W_t
\end{equation}
 where
$$A: [0,T]\times V\times \Omega\to V^*;\  \  B:
[0,T]\times V\times \Omega\to L_{2}(U;H)$$ are progressively
measurable. Suppose for a fixed $\alpha>1$  there exist constants
$\theta>0$, $K$ and a positive adapted process $f\in L^1([0,T]\times
\Omega; d
    t\times \mathbf{P})$ such that the
 following
 conditions hold for all $v,v_1,v_2\in V$ and $(t,\omega)\in [0,T]\times \Omega$.
\begin{enumerate}
    \item [$(A1)$] Hemicontinuity of $A$: The map
     $$ \lambda\mapsto { }_{V^*}\<A(t,v_1+\lambda
 v_2),v\>_V$$
     is
    continuous on $\mathbb{R}$.

    \item [$(A2)$] Monotonicity of $(A,B)$:
$$2{  }_{V^*}\<A(t,v_1)-A(t,v_2), v_1-v_2\>_V
    +\|B(t,v_1)-B(t,v_2)\|_{2}^2\\
     \le K\|v_1-v_2\|_H^2. $$

\item [$(A3)$] Coercivity of $(A,B)$:
    $$ 2{ }_{V^*}\<A(t,v), v\>_V +\|B(t,v)\|_{2}^2 +\theta
    \|v\|_V^{\alpha} \le f_t +K\|v\|_H^2.$$

\item[$(A4)$] Boundedness of $A$:
$$ \|A(t,v)\|_{V^*} \le f_t^{\alpha/(\alpha-1)} +
 K\|v\|_V^{\alpha-1}.$$
    \end{enumerate}
    Then for any $X_0\in L^2(\Omega\to H; \mathcal{F}_0;\mathbf{P})$,
    $(\ref{general equation})$
    has a unique solution $\{X_t\}_{t\in [0,T]}$ which is an adapted
 continuous process
    on $H$ such that $\mathbf{E}\int_0^T\|X_t\|_V^{\alpha}d t<\infty$
 and
$$\<X_t, v\>_H= \<X_0,v\>_H +\int_0^t{ }_{V^*}\<A(s,X_s),
    v\>_Vd s +\int_0^t \<B(s,X_s)d W_s, v\>_H$$
hold for all $v\in V$ and
$(t, \omega)\in [0,T]\times\Omega.$
 \end{lem}

Note that in order to using the coupling method, here we only
consider equation (\ref{1.1}) where the noise is the
additive type. We intend to establish Harnack inequality for
the associate transition semigroup
$$P_tF(x):= \E
F(X_t(x)),\ t\ge0, \ x\in H,$$
where $F$ is a bounded measurable function on $H$.
To define the intrinsic metric induced by $B_t$, we
need to assume  $B_t(\oo)$ is
non-degenerate for $t> 0$ and $\oo\in\OO$; that is, $B_t(\oo) y=0$
implies $y=0$. Then for $u\in V$
$$\|u\|_{B_t}:= \beg{cases} \|y\|_U, &\text{if}\ y\in U,\  B_t
y=u;\\
\infty, &\text{otherwise.}\end{cases}$$

\beg{thm}\label{T1.1} Assume $(A1)-(A4)$ hold for $(\ref{1.1})$ with
the coercivity exponent $\alpha$. Suppose
 there exist a constant $\sigma\ge2, \sigma>\alpha-2$
and continuous functions $\delta,\gamma,\xi\in
C[0,\infty)$ such that for any $t\geq0, \omega\in\Omega$ and $u,
v\in V$ we have
 \beq\label{c1}
2{ }_{V^*}\langle A(t,u)-A(t,v), u-v\rangle_V \leq -\delta_t
N(u-v)+\gamma_t \|u-v\|_{H}^2, \end{equation} \beq\label{c2}
N(u)\geq \xi_t \|u\|_{B_t}^\sigma\|u\|_{H}^{\alpha-\sigma},
 \end{equation}
where $N$
 is a positive real function on $V$ and $\xi,\delta$ are strictly positive on $[0,\infty)$,
 then $P_t$ is
 strong Feller operator for $t>0$, and for any $p>1$
and positive measurable function $F$ on $H$ we have
\beq\label{M}
 (P_tF(y))^p\leq
P_tF^p(x)\exp\[\frac{p}{p-1}C(t,\sigma)\|x-y\|_{H}^{2+\frac{2(2-\alpha)}{\sigma}}\],
\ x,y\in H,
\end{equation} where
$$C(t,\sigma)=\frac{2t^{\frac{\sigma-2}{\sigma}}(\sigma+2)^{2+\frac{2}{\sigma}}}{(\sigma+2-\alpha)^{2+\frac{2}{\sigma}}\[\int_0^t
(\delta_s\xi_s)^{\frac{1}{\sigma}}\exp(\frac{\alpha-2-\sigma
}{2\sigma}\int_0^s\gamma_udu)ds
\]^2} \  .$$
In particular, if  $\delta,\xi$ are time-independent, then
$$C(t,\sigma)=\frac{2(\sigma+2)^{2+\frac{2}{\sigma}}}{(\sigma+2-\alpha)^{2+\frac{2}{\sigma}}
(\delta\xi)^{\frac{2}{\sigma}}t^{\frac{\sigma+2}{\sigma}}}.$$
\end{thm}

\beg{rem} (1) Note that $(A1)-(A4)$  are assumed in  Theorem
\ref{T1.1} only for the existence and uniqueness of the strong
solution to (\ref{1.1}). One can replace those conditions by more
general ones in \cite{RRW,Zh}
 and obtain a  similar result.\\

(2) This theorem covers the main result in
\cite{W07} if we take $N(u)=\|u\|_V^{r+1}$ for  stochastic porous
media equations. Moreover, if we  take
$N(u)=\mathbf{m}(\mathbf{g}(u))$ for some Young function
$\mathbf{g}$, then this theorem can also be applied to stochastic
generalized porous media equations \cite{RRW}
in the framework of Orlicz space.\\

(3) This theorem can also be applied to many other types of
SPDE in \cite{PR,KR} which satisfy the strongly dissipative
 condition (\ref{c1}) (see section 3). For concrete examples in this paper we can consider
$N(u)=\|u\|_{V}^{\aa}$ for simplicity. In this case
(\ref{c1}) implies $(A2)$ and $(A3)$.
 Under (\ref{c1}) we have also
  established a stronger version of large deviation principle in \cite{L08}
for general SPDE with small multiplicative noise.\\

(4) Note  (\ref{c2}) implies that $V$ is contained in the range of
 $B_t$ (as a operator from $U$ to $H$) for fixed $t$ and $\omega$.
If we assume $N(u)=\|u\|_{V}^{\aa}$ and $V\equiv H$, then we know
 $B_t$ is a bijection map and
its inverse operator is also continuous from $H$ to $U$.
Since $B_t$ is a Hilbert-Schmidt operator, then
 $H$ and $U$ has to be finite dimensional space. In this case (\ref{c2})
holds provided $B_t$ are invertible.
\\

(5) The stochastic fast diffusion equations in \cite{RRW} does not
satisfy the assumption (\ref{c1}), but we have
 also
obtained the Harnack inequality, strong Feller property and heat
kernel estimates in \cite{LW}
 by using more delicate estimate. But we haven't obtained strong
contractive property (e.g.hyperboundedness) for the associated transition
semigroups
 in \cite{LW}
because of the weaker dissipativity of the drift.
\end{rem}

To apply Theorem \ref{T1.1} to obtain the heat kernel estimates,
ergodicity  and contractivity properties of $P_t$,
we only consider the deterministic and
time-homogenous case from now on. We first establish some
properties for  invariant measure.

\beg{thm}\label{T1}
 Suppose coefficients $A, B$ in $(\ref{1.1})$ are
deterministic and time-independent such that $(A1)$ and $(A4)$ hold.
Assume $(\ref{c1})$ hold
for $N(\cdot)=\|\cdot\|_V^\alpha$ and the embedding
$V\subseteq H$ is compact.

 (i) If $\gamma\leq0$ also holds in the case $\alpha\le2$,
then the Markov semigroup $\{P_t\}$ has an invariant probability measure
 $\m$, which satisfies
 $\mu\left(\|\cdot\|_V^\aa+e^{\varepsilon_0\|\cdot\|_H^\aa}\right)<\infty$
for some $\varepsilon_0>0$.

(ii) If $\alpha=2$ , then for any $x,y\in H$  we have
$$\|X_t(x)-X_t(y)\|_H^2\le e^{(\gamma-c_0\delta)t}\|x-y\|_H^2, \ t\ge 0,$$
where $c_0$ is the constant such that $\|\cdot\|_V^2\ge c_0\|\cdot\|_H^2$ hold.

Moreover, if $\gamma<c_0\delta$, then there exists a unique
  invariant measure $\mu$ of $\{P_t\}$  and
for any Lipschitz continuous function $F$ on $H$  we have
\beq\label{exponentially decay}
|P_tF(x)-\mu(F)|\le  \text{Lip}(F) e^{-(c_0\delta-\gamma)t/2}(\|x\|_H+C),
\   x\in H,
\end{equation}
 where $C>0$ is a constant and $\text{Lip}(F)$ is the Lipschitz constant of $F$.

(iii) If $\alpha>2$ and $\gamma\leq0$, then there exists a constant $C$ such that
$$ \|X_t(x)-X_t(y)\|_H^2\le \|x-y\|_H^2 \wedge \left\{ C t^{-\frac{2}{\alpha-2}}\right\},\ t>0,\ x,y\in H $$
where $X_t(y)$ is the solution to $(\ref{1.1})$ with the starting point $y$.

Therefore,
 $\{P_t\}$ has a unique  invariant measure
$\mu$ and
for any Lipschitz continuous function $F$ on $H$ we have
\beq\label{algebraically decay}
\sup_{x\in H}|P_tF(x)-\mu(F)|\le C \text{Lip}(F) t^{-\frac{1}{\alpha-2}},
\ t>0.
\end{equation}
\end{thm}

\begin{rem}
 (\ref{algebraically decay}) describes the algebraically
convergence rate of the transition semigroup to the invariant measure.
In particular, if $B=0$ and Dirac measure at $0$ is the unique
invariant measure of  $\{P_t\}$, then we can take $F(x)=\|x\|_H$
in (\ref{algebraically decay}) and have
$$\sup_{x\in H} \|X_t(x)\|_H\le C t^{-\frac{1}{\alpha-2}}, \
t>0. $$
Hence it gives the decay estimate of the solution
to a large class of deterministic evolution equations.
These results coincide with some well-known decay estimates
 in PDE theory, e.g. the optimal decay of the solution
to the
 classical porous medium equation in \cite{AP,DRRW}.
We refer to section 3 for more examples.
\end{rem}

 We  recall that  $\{P_t\}$ is called
(topologically)
 irreducible
if $P_t1_M(\cdot)>0$ on $H$ for any $t>0$ and
nonempty open set $M$. If $\{P_t\}$ is a semigroup
 defined on $L^2(\mu)$,
then $\{P_t\}$ is called hyperbounded semigroup if
 $\|P_t\|_{L^2(\m)\rightarrow L^4(\m)}<\infty$ for some $t>0$;
$\{P_t\}$ is called ultrabounded semigroup if
 $\|P_t\|_{L^2(\m)\rightarrow L^\infty(\m)}<\infty$ for any $t>0$.

 \beg{thm}\label{T1.2}
 Suppose coefficients $A, B$ in $(\ref{1.1})$ are
deterministic and time-independent such that all assumptions
in Theorem $\ref{T1.1}$ hold for $N(\cdot)=\|\cdot\|_V^\alpha$.

(i) $\{P_t\}$ is  irreducible and has a
unique  invariant measure $\mu$
with full support on $H$.
 Moreover, $\mu$ is strong mixing and
for any probability measure $\nu$ on $H$ we have
$$\lim_{t\rightarrow\infty}\parallel P_t^*\nu-\mu\parallel_{var}=0 ,$$
where $\parallel\cdot\parallel_{var}$ is the total variation norm
 and $P_t^*$ is the adjoint operator of
$P_t$.

 (ii) For any $x\in H$,
$t>0$ and  $p>1$, the transition density
$p_t(x,y)$ of $P_t$ w.r.t $\m$ satisfies
$$\|p_t(x,\cdot)\|_{L^p(\m)}\leq\left\{\int_H
\exp\left[-pC(t,\sigma)\|x-y\|_{H}^{2+\frac{2(2-\alpha)}{\sigma}}\right]\mu
(dy)\right\}^{-\frac{p-1}{p}}.$$

(iii) If $\aa=2$, then $P_t$ is hyperbounded
and compact on
$L^2(\m)$ for some $t>0$.

(iv) If $\aa>2$,
then $P_t$ is ultrabounded and
compact on $L^2(\m)$ for any $t>0$. Moreover, there exists a
constant $C>0$ such that
$$\|P_t\|_{L^2(\m)\rightarrow L^\infty(\m)}
 \leq\exp\left[C(1+t^{-\frac{\aa}{\aa-2}})\right],
\ t>0.$$
\end{thm}

\begin{rem}
Based on the Harnack inequality, the  irreducibility
is derived very easily for the transition semigroup. Then
according to  Doob's theorem (cf.\cite{MS,H03}) one
can derive the
 uniqueness of
invariant measures and some ergodic properties for the associated transition
 semigroups. Comparing with the uniqueness result for  invariant measure
in Theorem \ref{T1}, we do not need to
assume $\gamma\le 0$ or $\gamma<c_0\delta$
in this case.

\end{rem}

Let $\L_p$ be the generator of the semigroup $\{P_t\}$ in $L^p(\mu)$. We say that $\L_p$ has the
spectral gap in $L^p(\mu)$ if there exists $\gamma>0$ such that
$$ \sigma(\L_p)\cap\{\lambda: Re \lambda>-\gamma\}=\{0\} $$
where $\sigma(\L_p)$ is the spectrum of $\L_p$. The largest constant $\gamma$ with this property is denoted by $gap(\L_p)$.

\begin{thm}\label{T1.3}
 Suppose  all assumptions in Theorem \ref{T1.2} hold and $\mu$ denotes the unique invariant measure of $\{P_t\}$.

(i) If $\alpha=2$ and $\gamma<c_0\delta$, then the Markov semigroup $\{P_t\}$ is $V$-uniformly ergodic, i.e. there exist $C,\eta>0$
such that for all $t\ge0$ and $x\in H$
$$\sup_{\|F\|_V\le 1}|P_tF(x)-\mu(F)|\le C V(x)e^{-\eta t},  $$
where   we can take $V(x)=1+\|x\|_H^2$ and $V(x)=e^{\varepsilon_0\|x\|_H^2}$ for some  small constant $\varepsilon_0>0$,
$$
 \|F\|_V:=\sup_{x\in H}\frac{|F(x)|}{V(x)}<\infty.$$
Moreover, if $P_t$ is symmetric on $L^2(\mu)$ for all $t\ge 0$, then  we have
$$  \|P_tF-\mu(F)\|_{L^2(\mu)}\le  e^{-\eta t}\|F\|_{L^2(\mu)},\ F\in L^2(\mu),\ t\ge 0.    $$

(ii) If $\alpha>2$, then the Markov semigroup $\{P_t\}$ is uniformly exponential ergodic, i.e. there exist $C,\eta>0$
such that for all $t\ge0$ and $x\in H$
$$\sup_{\|F\|_\infty\le 1}|P_tF(x)-\mu(F)|\le C e^{-\eta t}.  $$
Moreover, for each $p\in(1,\infty]$ we have
$$  \|P_tF-\mu(F)\|_{L^p(\mu)}\le C_p e^{-(p-1)\eta t/p}\|F\|_{L^p(\mu)},\  F\in L^p(\mu),\  t\ge 0,    $$
and
$$ gap(\L_p)\ge \frac{(p-1)\eta}{p},$$
where $C_p$ is a constant and we set $\frac{p-1}{p}=1$ if $p=\infty$ by convention.
\end{thm}

\begin{rem} The $V$-uniformly ergodicity implies that for any probability measure
$\nu$ on $H$ we have
\begin{equation*}\begin{split}
\|P_t^*\nu-\mu\|_{var}&\le \int_H \|P(t,x,\cdot)-\mu\|_{var} \nu(\d x)\\
&\le \int_H \sup_{\|\varphi\|_V\le 1} \left|P_t\varphi(x)-\mu(\varphi)  \right| \nu(\d x)\\
&\le \int_H CV(x)e^{-\eta t} \nu(\d x)
= C \nu(V) e^{-\eta t},\ t\ge0.
\end{split}\end{equation*}
 And it is easy to show that the uniformly exponential ergodicity is equivalent to
$$ \|P_t^*\nu-\mu\|_{var}\le Ce^{-\eta t},\ t\ge0, $$
since we have
$$\sup_{\nu} \|P_t^*\nu-\mu\|_{var}=\sup_{x\in H}  \|P_t(x,\cdot)-\mu\|_{var}
=\frac{1}{2}\sup_{\|f\|_\infty\le1}\|P_tf-\mu(f)\|_{\infty}.  $$
\end{rem}

The paper is organized as follows.
The main theorems are proved in  section 2.  To apply
the main results, one has to verify condition (\ref{c1})
and (\ref{c2}). For this purpose a crucial inequality is proved as a
lemma in section 3. Then some concrete examples are discussed  as applications.

\section{Proofs of the main theorems}
\subsection{Proof of Theorem \ref{T1.1}}
The main
techniques in the proof are a coupling argument and Girsanov transformation
 in infinite
dimensional space (cf.\cite{LW,W07}). The coupling method dates back to D\"{o}blin's work
\cite{D38} on Markov chains
and it is one of the main tools in particle systems (cf.\cite{C04}).
The first use of coupling for SPDE
 up to our knowledge is due to Mueller \cite{M93}, who used
this technique to prove
 the uniqueness of invariant measures for the stochastic
heat equation. We refer to some review papers
 \cite{M03,MS99,H03} on this subject
for more references.

The coupling we used here,
which only depends on the natural distance between two marginal processes,
is a modification of the argument in \cite{ATW}.
Such a stronger
 Harnack
inequality (the estimate only depending on the usual norm) will provide more
information such as the strong Feller property and the hyper- or
ultrabounded property of the transition semigroups.
In order to make the proof easier to understand, we first
describe the main ideas and steps.

To prove the Harnack inequality for the transition semigroup $\{P_t\}$,
it suffices to construct
a coupling $(X_t, Y_t)$, which is a continuous adapted process on
$H\times H$ such that
\\
\noindent(i) $X_t$ solves (\ref{1.1}) with $X_0=x$;\\
(ii) $Y_t$ solves the following equation
$$\d Y_t = A(t,Y_t)\d t +B_t \d \tilde{W}_t,\ Y_0=y$$ for another
 cylindrical Brownian motion
$\tilde{W}_t$ on $U$ under a weighted probability measure
$R\mathbf{P}$, where $\tilde{W}_t$ as well as the density $R$ will
be constructed  by a
Girsanov transformation;\\
(iii) $X_T=Y_T, a.s.$
\\

 As soon as (i)-(iii) are satisfied, then we have
\beq\label{h}\begin{split}
P_TF(y)&=\mathbf{E}RF(Y_T)=\mathbf{E}RF(X_T)\\
&\leq(\mathbf{E}R^{p/(p-1)})^{(p-1)/p}(\mathbf{E}F^p(X_T))^{1/p}\\
       &=(\mathbf{E}R^{p/(p-1)})^{(p-1)/p}(P_TF^p(x))^{1/p},
\end{split}\end{equation}
which implies the desired Harnack inequality provided
$\mathbf{E}R^{p/(p-1)}<\infty$.

Now we constract the coupling process $Y_t$.  We first take
$\varepsilon\in(0, 1), \ \beta\in \mathbf{C}([0,\infty);
\mathbb{R}_+)$ and consider the equation
 \beq\label{2.1} \d
Y_t=\left(A(t,Y_t)+\frac{\beta_t(X_t-Y_t)}{\xy_H
^\varepsilon}\mathbf{1}_{\{t<\tau\}}\right)\d t+B_t\d W_t, \ Y_0=y,
\end{equation}
where $X_t:=X_t(x)$ and $\tau:=\inf\{t\geq0: X_t=Y_t\}$ is the
coupling time.

 According to Lemma \ref{L0} we can prove
that (\ref{2.1}) also has a unique strong solution $Y_t(y)$ by using
 a
similar argument in \cite[Theorem A.2]{W07} (in fact, one can prove
the added drift is also monotone). Then by (\ref{c1}) we have
$$ \xy_H^2
\leq
\|X_s-Y_s\|_H^2+\int_s^t\left(-\delta_uN(X_u-Y_u)+\gamma_u\|X_u-Y_u\|_H^2
-\beta_u\|X_u-Y_u\|_H^{2-\varepsilon}\mathbf{1}_{\{u<\tau\}}\right)\d
 u$$
for all $0\leq s\leq t$.
 Hence we have $X_t=Y_t$ for $t\geq
\tau$ by using  Gronwall's lemma.

And it is  easy to show that
\beq\label{2.0} e^{-\int_0^t\gamma_s\d s}\xy_H^2\leq
\|x-y\|_H^2-\int_0^t e^{-\int_0^u\gamma_sds}
\(\delta_uN(X_u-Y_u)+\beta_u\|X_u-Y_u\|_H^{2-\varepsilon}
\mathbf{1}_{\{u<\tau\}}\)\d u.
\end{equation}
First, we will prove  the coupling time $\tau\leq T \
a.s.$ by choosing  $\beta_t$ appropriately   in (\ref{2.1}).

\beg{lem}\label{L1}
 If $\beta$ satisfies
 $\int_0^T\beta_t e^{-\frac{\varepsilon}{2}\int_0^t\gamma_s\d s}\d
 t\geq\frac{2}{\varepsilon}\|x-y\|_H^\varepsilon$,
 then $X_T=Y_T ,\
 a.s.$
 \end{lem}

 \beg{proof} By (\ref{2.0}) and  the chain rule we have
$$\left\{e^{-\int_0^t\gamma_s\d s} \xy_H^2\right\}^{\varepsilon/2}
 \leq \|x-y\|_H^\varepsilon-\frac{\varepsilon}{2} \int_0^t\beta_s
e^{-\frac{\varepsilon}{2}\int_0^s\gamma_u\d u}\d s ,\
t\leq\tau\wedge
 T .$$
If $T<\tau(\omega_0)$ for some $\omega_0\in \Omega$, then by
taking $t=T$ and using the assumption we have
$$e^{-\frac{\varepsilon}{2}\int_0^T\gamma_s\d s}
\|X_T(\omega_0)-Y_T(\omega_0)\|_H^\varepsilon \leq
\|x-y\|_H^\varepsilon-\frac{\varepsilon}{2}\int_0^T\beta_t
e^{-\frac{\varepsilon}{2}\int_0^t\gamma_s\d s}\d t\leq 0.$$ This
implies
 $X_T(\omega_0)=Y_T(\omega_0)$, which
contradicts with the assumption $T<\tau(\omega_0)$.

 Hence
$\tau\leqslant T,\  a.s.$ The proof is complete.
\end{proof}

\noindent $\mathbf{Proof\ of\ Theorem\ \ref{T1.1}}:$\ Let
$\varepsilon=1-\frac{\alpha}{\sigma+2}\in(0, 1)$, then by (\ref{2.0}) and
(\ref{c2}) we have \beq\beg{split}\label{2.2}
 \d\left\{\xy_H^2 e^{-\int_0^t\gamma_s\d s}
 \right\}^\varepsilon&\leq-\varepsilon\delta_t
 e^{-\varepsilon\int_0^t\gamma_s\d
s}\xy_H^{2(\varepsilon-1)}N(X_t-Y_t)
 \d t\\
&\leq-\varepsilon\delta_t\xi_t  e^{-\varepsilon\int_0^t\gamma_s\d s}
\frac{\xy_{B_t}^\sigma}{\xy_H^{2+\sigma-\alpha-2\varepsilon}}\d t\\
&=-\varepsilon\delta_t\xi_t  e^{-\varepsilon\int_0^t\gamma_s\d
 s}\frac{\xy_{B_t}^\sigma}{\xy_H^{\sigma\varepsilon}}\d t\\
&=-\frac{\beta_t^\sigma\xy_{B_t}^\sigma}{c^\sigma\xy_H^{\sigma\varepsilon}}\d
t,
\end{split}\end{equation}
 where
$$\beta_t^\sigma=c^\sigma\varepsilon\delta_t\xi_t e^{-\varepsilon\int_0^t\gamma_s\d
s},\ \
c=\frac{2\|x-y\|_H^\varepsilon}{\varepsilon\int_0^T(\varepsilon\delta_t\xi_t)^{\frac{1}{\sigma}}
 e^{-(\frac{1}{2}+\frac{1}{\sigma})\varepsilon\int_0^t\gamma_s\d s}\d
 t}.$$
Let$$\zeta_t:=\frac{\beta_tB_t^{-1}(X_t-Y_t)}{\xy_H^\varepsilon}\mathbf{1}_{\{t<\tau\}}.$$
 By using
H\"{o}lder's inequality and (\ref{2.2}) we obtain \beq\beg{split}
\int_0^T\|\zeta_t\|_U^2\d
 t&=\int_0^T\frac{\beta_t^2\xy_{B_t}^2}{\xy_H^{2\varepsilon}}\d t\\
&\leq
 T^{\frac{\sigma-2}{\sigma}}\(\int_0^T\frac{\beta_t^\sigma\xy_{B_t}^\sigma}{\xy_H^{\sigma\varepsilon}}\d t\)
^{\frac{2}{\sigma}}\\
&\leq
T^{\frac{\sigma-2}{\sigma}}\(c^\sigma\|x-y\|_H^{2\varepsilon}\)^{\frac{2}{\sigma}}.
\end{split}\end{equation}
Hence we have
 \beq\beg{split}\label{2.3}
 \mathbf{E}\exp\[\frac{1}{2}\int_0^T\|\zeta_t\|^2_U\d t\]<\infty.
\end{split}\end{equation}
Therefore, we can rewrite
(\ref{2.1}) as
$$\d Y_t=A(t,Y_t)\d t+B_t\d\tilde{W}_t,\  Y_0=y$$
where
$$\tilde{W}_t:=W_t+\int_0^t\zeta_s\d s.$$

By (\ref{2.3}) and the Girsanov theorem (e.g.\cite[Th
 10.14,\ Prop.10.17]{DaZa}) we know that
$\{\tilde{W}_t\}$ is a cylindrical Brownian motion on
$U$ under the weighted probability measure $R\mathbf{P}$, where
$$R:=\exp\[\int_0^T\langle\zeta_t, \d
 W_t\rangle-\frac{1}{2}\int_0^T\|\zeta_t\|_U^2\d t\] .$$
Therefore, the distribution of $\{Y_t(y)\}_{t\in[0,T]}$ under
$R\mathbf{P}$ is  same with the distribution of
$\{X_t(y)\}_{t\in[0,T]}$ under $\mathbf{P}$.

Let $p^\prime=\frac{p}{p-1}$, then for any $q>1$
 \beq\label{estimate  of R}\beg{split}
\mathbf{E}R^{p^\prime}&=\exp\[p^\prime\int_0^T\langle\zeta_t,
\d W_t\rangle-\frac{p^\prime}{2}\int_0^T\|\zeta_t\|_U^2\d t\]\\
&\leq \[\mathbf{E}\exp( qp^\prime\int_0^T\langle \zeta_t, \d
W_t\rangle-\frac{q^2(p^\prime)^2}{2}\int_0^T\|\zeta_t\|_U^2\d t)
\]^{\frac{1}{q}}\\
& \ \ \ \cdot
\[\mathbf{E}\exp(\frac{qp^\prime(qp^\prime-1)}{2(q-1)}\int_0^T\|\zeta_t\|_U^2\d
 t)\]^{\frac{q-1}{q}}\\
&\leq
 \[\mathbf{E}\exp(\frac{qp^\prime(qp^\prime-1)}{2(q-1)}\int_0^T\|\zeta_t\|_U^2\d t)\]^{\frac{q-1}{q}}\\
&\leq
\exp\[\frac{p^\prime(qp^\prime-1)}{2}T^{\frac{\sigma-2}{\sigma}}
\(c^\sigma\|x-y\|_H^{2\varepsilon}\)^{\frac{2}{\sigma}}\].
\end{split}\end{equation}
 By taking $q\downarrow1$ we have
\beq\beg{split}
(P_T(y))^p&\leq P_TF^p(x)(\mathbf{E}R^{p\prime})^{p\prime-1}\\
&\leq
P_TF^p(x)\exp\[\frac{p}{p-1}C(t,\sigma)\|x-y\|_{H}^{2+\frac{2(2-\alpha)}{\sigma}}\],
\end{split}\end{equation}
where
$$C(t,\sigma)=\frac{2t^{\frac{\sigma-2}{\sigma}}(\sigma+2)^{2+\frac{2}{\sigma}}}{(\sigma+2-\alpha)^{2+\frac{2}{\sigma}}\[\int_0^t
(\delta_s\xi_s)^{\frac{1}{\sigma}}\exp(\frac{\alpha-2-\sigma
}{2\sigma}\int_0^s\gamma_u\d u)\d s
\]^2} \  .$$

From (\ref{estimate  of R}) we know that $R$ is uniformly
integrable,  then by the dominated
convergence theorem we have
$$\lim_{y\to x} \mathbf{E}|R-1| =\mathbf{E}\lim_{y\to x}|R-1|=0.$$
Hence
$$ |P_TF(y)-P_TF(x)| =|\mathbf{E}RF(X_T)-\mathbf{E}F(X_T)| \le \|F\|_\infty
\mathbf{E}|R-1|\rightarrow 0(y\rightarrow x).$$
This implies $P_TF\in
C_b(H)$. Therefore,  $P_T$
is strong Feller operator. \qed

\subsection{Proof of Theorem \ref{T1}}

 $(i)$ In the present case, $\{P_t\}$ is a Markov semigroup (cf.\cite{KR,PR}).
 The existence  of an invariant measure $\mu$ can be proved by
 the standard Krylov-Bogoliubov procedure (cf.\cite{PR,W07}).
 Let
 $$\mu_n:=\frac{1}{n}\int_0^n\delta_0P_tdt,\  n\geq1 ,$$
 where $\delta_0$ is the Dirac measure at $0$.
Recall $X_t(y)$ is the solution to (\ref{1.1}) with the starting point $y$, then by
$(\ref{c1})$ and the Gronwall Lemma
$$\|X_t(x)-X_t(y)\|_H^2\le e^{\gamma t}\|x-y\|_H^2, \ \  \forall x,y\in H.$$
This implies that $P_t$ is a  Feller semigroup.

Hence for the existence of an invariant measure, it is well-known that one
 only needs to verify the tightness of
 $\{\mu_n :  n\geq1\}$.

Since  $\gamma\le 0$ in the case $\alpha\le 2$, then by (\ref{c1}) and $(A4)$
  we have
 \beq\label{inequality 1}\beg{split}
  2{ }_{V^*}<A(x), x>_V&\leq -\delta\|x\|_{V}^{\alpha}+ \gamma\|x\|_H^2+
2{~}_{V^*}\langle
 A(0),x\rangle_V\\
&\leq \theta_2-\theta_1\|x\|_{V}^\alpha
\end{split}\end{equation}
holds for some constant $\theta_1, \theta_2>0$.
By using the It\^{o} formula  we have
\beq\label{estimate for invariant measure}
\x_H^2
\leq \|x\|_H^2+\int_0^t (c-\theta_1\|X_s\|_{V}^{\alpha})ds+
2\int_0^t\langle X_s, BdW_s\rangle_H,
 \end{equation}
  where $c>0$ is some constant which may change from
line to line.

Note that $ M_t:=\int_0^t\< X_s, BdW_s\>_H$ is a
 martingale, then  (\ref{estimate for invariant measure})
implies that
\beq\label{estimate of mu_n}
\mu_n(\|\cdot\|_V^\alpha)=\frac{1}{n}
\int_0^n\mathbf{E}\|X_t(0)\|_{V}^{\alpha}dt\leq\frac{c}{\theta_1},
\  n\geq1.
\end{equation}
Since the embedding  $V \subseteq H$
is compact, then
 for any constant $K$ the set
$\{x\in H: \ \|x\|_{V}\le K\}$
is relatively compact
in $H$.  Therefore, (\ref{estimate of mu_n}) implies that
$\{\mu_n\}$ is tight, hence the limit of a convergent subsequence
provides an invariant measure $\mu$ of $\{P_t\}$.

Now  we need to prove the concentration property of $\mu$.
 If  $\varepsilon_0$ is small enough,
then by (\ref{estimate for invariant measure})
 and  It\^{o}'s formula
\begin{equation}\begin{split}
\label{exponential estimate}
e^{\varepsilon_0\x_H^\alpha}\leq&  e^{\varepsilon_0\|x\|_H^\alpha}
+\int_0^t\left(c
-\theta_1\|X_s\|_V^\alpha+\alpha\varepsilon_0\|B\|_2^2\|X_s\|_H^\alpha\right)\frac{\alpha\varepsilon_0}{2}
\|X_s\|_H^{\alpha-2}e^{\varepsilon_0\|X_s\|_H^\alpha}\d s\\
&+ \alpha\varepsilon_0\int_0^t
\|X_s\|_H^{\alpha-2}e^{\varepsilon_0\|X_s\|_H^\alpha}\langle X_s, BdW_s \rangle_H \\
\leq&  e^{\varepsilon_0\|x\|_H^\alpha}
+\int_0^t\left(c
-c_1\|X_s\|_H^\alpha\right)\frac{\alpha\varepsilon_0}{2}
\|X_s\|_H^{\alpha-2}e^{\varepsilon_0\|X_s\|_H^\alpha}\d s\\
& + \alpha\varepsilon_0\int_0^t
\|X_s\|_H^{\alpha-2}e^{\varepsilon_0\|X_s\|_H^q}\langle X_s, BdW_s \rangle_H\\
\leq &  e^{\varepsilon_0\|x\|_H^\alpha}
+\int_0^t\left(c_2
-c_3 e^{\varepsilon_0\|X_s\|_H^\alpha}\right)\d s
+ \alpha\varepsilon_0\int_0^t
\|X_s\|_H^{\alpha-2}e^{\varepsilon_0\|X_s\|_H^\alpha}\langle X_s, BdW_s \rangle_H
\end{split}
\end{equation}
holds for some positive
constants $c, c_1, c_2$ and $c_3$.
Therefore
 $$\mu_n(e^{\varepsilon_0\|\cdot\|_H^\alpha})=
\frac{1}{n}\int_0^n\mathbf{E}e^{\varepsilon_0\|X_t(0)\|_H^\alpha}\d t
 \leq \frac{1}{c_3n}+ \frac{c_2}{c_3},\ n\ge 1.$$
 Hence we have $\mu(e^{\varepsilon_0\|\cdot\|_H^\alpha})<\infty$ for some $\varepsilon_0>0$.
In particular, this implies  $\mu(\|\cdot\|_H^2)<\infty$.

By (\ref{estimate for invariant measure}) there also exists a constant $C$ such that
$$\mathbf{E} \int_0^1\|X_t(x)\|_{V}^\alpha dt \le C(1+\|x\|_H^2),\  \forall x\in H. $$
Therefore
$$\mu(\|\cdot\|_{V}^\alpha)=\int_H\mu(dx)\int_0^1\mathbf{E}(\|X_t(x)\|_{V}^\alpha) dt
\le C+C\int_H\|x\|_H^2\mu(dx)<\infty.$$




$(ii)$ If $\alpha=2$, then for any  $x,y\in H$
$$\|X_t(x)-X_t(y)\|_H^2\le \|x-y\|_H^2
+\int_0^t\left(-\delta\|X_s(x)-X_s(y)\|_V^2+\gamma
\|X_s(x)-X_s(y)\|_H^2\right)d s.$$
By the Gronwall lemma we have
$$ \|X_t(x)-X_t(y)\|_H^2\le e^{(\gamma-c_0\delta)t} \|x-y\|_H^2 , \
\forall  x,y\in H. $$
If $\gamma<c_0\delta$, then (\ref{inequality 1}) still holds. Hence
$\{P_t\}$ has an invariant measure by repeating the argument in $(i)$.
And we also have
$$ \lim_{t\rightarrow \infty} \|X_t(x)-X_t(y)\|_H=0, \
\forall  x,y\in H.    $$
By the dominated convergence theorem we know
 for any invariant measure $\mu$ and for any bounded continuous function $F$
$$|P_tF(x)-\mu(F)|\le \int_{H} \mathbf{E}
|F(X_t(x))-F(X_t(y))|\mu(dy)\rightarrow 0 (t\rightarrow\infty).$$
This implies the uniqueness of invariant measures.

We denote the invariant measure by $\mu$. By $(i)$ we know $\mu(\|\cdot\|_H^2)<\infty$,
hence for any bounded Lipschitz function $F$ on $H$ we have
$$\begin{aligned}
|P_tF(x)-\mu(F)|& \le \int_{H} \mathbf{E}
|F(X_t(x))-F(X_t(y))|\mu(dy)\\
& \le \text{Lip}(F)e^{(\gamma-c_0\delta)t/2}
 \int_H\|x-y\|_H\mu(dy)\\
&\le \text{Lip}(F)e^{(\gamma-c_0\delta)t/2}\left(\|x\|_H+C\right), \ x\in H,
 \end{aligned}$$
where  $C>0$ is a constant.

$(iii)$ If $\alpha>2$ and $\gamma\le0$, then there
exists a constant $c>0$  such that
$$\|X_t(x)-X_t(y)\|_H^2\le \|x-y\|_H^2
-c\int_0^t\|X_s(x)-X_s(y)\|_H^{\alpha} ds, \
t\geq 0.$$
Suppose $h_t$ solves the equation
\beq\label{compare function}
h_t^\prime=-ch_t^{\frac{\alpha}{2}},\
 h_0=(\|x-y\|_H +\varepsilon)^2,
\end{equation}
where $\varepsilon$ is a positive constant.
Then by a standard comparison argument  we have
\begin{equation}\label{comparison estimate}
 \|X_t(x)-X_t(y)\|_H^2\le h_t\le Ct^{-\frac{2}{\alpha-2}},
\end{equation}
where $C>0$ is a constant.
In fact, we can define
 $$\varphi_t:=h_t-\|X_t(x)-X_t(y)\|_H^2, \ \
\tau:=\inf\{t\ge0:\ \varphi_t<0 \}.$$
If $\tau<\infty$, then we know $\varphi_\tau\le 0$
by the continuity.

By the mean-value theorem
we have
$$\aligned
\varphi_t &\ge \varphi_0-c\int_0^t
\left(h_s^{\frac{\alpha}{2}}-\|X_s(x)-X_s(y)\|_H^\alpha\right)d s\\
&\ge \varepsilon^2-K\int_0^t\varphi_s d s, \ 0\le t\le\tau,
\endaligned$$
where $K>0$ is some constant. Then by the Gronwall lemma we have
$$\varphi_\tau\ge \varepsilon^2 e^{-K\tau}>0,$$
which is contradict to $\varphi_\tau\le 0$.
Hence (\ref{comparison estimate}) holds.

Therefore, for any $x\in H$ and bounded Lipschitz function $F$ on $H$, we have
$$|P_tF(x)-\mu(F)|\le \int_H \textbf{E}
|F(X_t(x)-F(X_t(y)))|\mu(dy)\le C\text{Lip}(F)t^{-\frac{1}{\alpha-2}}.$$
Hence (\ref{algebraically decay}) holds and the uniqueness of
invariant measures also follows.
\qed

\subsection{Proof of Theorem \ref{T1.2}}

$(i)$ By the definition of $\|\cdot\|_B$ and (\ref{c2}), for any constant $K$
there exists $\overline{K}>0$ such that
$$\begin{aligned}
\{x\in H: \ \|x\|_B\le K\}& \subseteq \{Bu:\ u\in U; \  \|u\|_U\le \overline{K}\};\\
\{x\in H:\  \|x\|_V\leq K\}& \subseteq
\{x\in H: \|x\|_B\leq
\overline{K}\}.
\end{aligned}$$
Since $B$ is a Hilbert-Schmidt (hence compact) operator,
then   the following set
$$\{x\in
H:\ \|x\|_V\leq K\}$$
is relatively compact in $H$ for any constant $K$,
$i.e.$ the embedding  $V \subseteq H$ is compact.
Hence
$\{P_t\}$ has an invariant measure according to Theorem \ref{T1}.

Suppose $\mu$ is an invariant
measure of $P_t$, then by taking $p=2$ in (\ref{M}) we have
\beq\beg{split}\label{density}
 & (P_t 1_M(x))^2 \int_H
e^{-2C(t,\sigma)\|x-y\|_H^{2+\frac{2(2-\alpha)}{\sigma}}}\mu(d
y) \\
&\le \int_HP_t 1_M(y)\mu(dy)= \mu(M),
\end{split}\end{equation}
 where $M$ is a Borel set on $H$.
 Hence the transition kernel $P_t(x,dy)$ is absolutely continuous
w.r.t. $\mu$, and we denote the density by  $p_t(x,y)$.

If $\mu$ does not have full support on $H$,
this means there exist $x_0\in H$ and $r>0$ such
that
$$B(x_0;r):= \{y\in H: \|y-x_0\|_H\le r\}$$
 is a null set of
$\mu$. Then (\ref{density}) implies that
 $P_t(x_0, B(x_0;r))=0$, $i.e.$
$$\mathbf{P}\left(X_t(x_0)\in B(x_0;r)\right)=0,\ \ \ t>0.$$
Since $X_t(x_0)$ is a continuous process on $H$,  we have
$\mathbf{P}\left(X_0\in B(x_0;r)\right)=0$, which is
contradict with $X_0=x_0$.

Therefore,
$\mu$ has full support on $H$.

According to the Harnack inequality (\ref{M}) we have
$$(P_t1_M)^p(x_0)\leq
P_t1_M(x)\exp\[\frac{p}{p-1}C(t,\sigma)\|x-x_0\|_{H}^{2+\frac{2(2-\alpha)}{\sigma}}\]
,\ x,x_0\in H.$$
Therefore,  to prove the irreducibility, one only has to show
for any given nonempty
 open set $M$ and $t>0$,
 there exists $x_0\in H$ such that $P_t1_M(x_0)>0$ .

Note that the
 full support property of $\mu$ implies
 $$\int_H P_t1_M(x)\mu(\d x)=\int_H 1_M(x)\mu(\d x)=\mu(M)>0.$$
So $P_t1_M(\cdot)$ cannot be the zero function. Therefore $\{P_t\}$
is irreducible.

Since $\{P_t\}$ have also the strong Feller property, then the uniqueness of invariant
measures follows from the classical Doob theorem \cite{Do} (or see
\cite[Theorem 2.1]{H03}).

 Note that the solution
 has continuous paths on $H$, then the other assertions follow from the
 general result in the ergodic theory (cf.\cite[Theorem 2.2 and Proposition 2.5]{Se}, \cite{MS}).

$(ii)$ For any $p>1$ and  nonnegative measurable function $f$ with
$\mu(f^{p/(p-1)})\le 1$,   by replacing $p$ with $p/(p-1)$ in (\ref{M})
we have
$$ \big(P_t f(x)\big)^{p/(p-1)}\le \big(P_t
f^{p/(p-1)}(y)\big)\exp\big[pC(t,\sigma)
\|x-y\|_H^{2+\frac{2(2-\alpha)}{\sigma}}\big],\ \ \ x,y\in H.$$
Taking integration w.r.t. $\mu(dy)$ on both sides we have
$$\big(P_t f(x)\big)^{p/(p-1)}\int_H e^{-pC(t,\sigma)
\|x-y\|_H^{2+\frac{2(2-\alpha)}{\sigma}}}\mu(d y)\le
\mu(f^{p/(p-1)})\le 1.$$
This implies that
$$P_t f(x) \le \bigg(\int_H
e^{-pC(t,\sigma) \|x-y\|_H^{2+\frac{2(2-\alpha)}{\sigma}}}\mu(d y)
\bigg)^{-(p-1)/p}.$$
Note that
$$P_tf(x)=\int_H f(y)P_t(x,dy)=\int_H f(y)p_t(x,y)\mu(dy),$$
hence we have
\begin{equation*}\begin{split}
\|p_t(x,\cdot)\|_{L^p(\mu)}&=
\sup_{\|f\|_{L^q(\mu)}\le 1}\left| \int_H f(y)p_t(x,y)\mu(dy)\right|\\
& \le \bigg(\int_H
e^{-pC(t,\sigma) \|x-y\|_H^{2+\frac{2(2-\alpha)}{\sigma}}}\mu(d y)
\bigg)^{-(p-1)/p},
\end{split}\end{equation*}
where $q=p/(p-1)$.

$(iii)$
By  (\ref{M})  there exists a constant $c>0$
such that \beq
 (P_t
f)^2(x)\exp\bigg[-\frac{c\|x-y\|_H^{2+\frac{2(2-\alpha)}{\sigma}}}{t^{\frac{\sigma+2}{\sigma}}}\bigg]\le
P_t
 f^2(y),\ \ \ x,y\in H,\ t>0.
\end{equation}
Integrating on both sides w.r.t. $\mu(d y)$,
for $f\in L^2(\mu)$ with $\mu(f^2)=1$ we have
\beq\beg{split}\label{semigroup inequality}
  (P_t f)^2(x)
\le \frac{1}{\mu(B(0,1))} \exp\Big[
\frac{c(\|x\|_H+1)^{2+\frac{2(2-\alpha)}{\sigma}}}{t^{\frac{\sigma+2}{\sigma}}}\Big],\
\ \ x\in H, t>0, \end{split}\end{equation}
 where $B(0;1):= \{y\in H: \|y\|_H\le 1\}$ and
$\mu\left(B(0;1)\right)>0$.

 If $\aa=2$,  then there exists $C>0$ such that
$$\int_H (P_tf)^4(x)\mu(d x) \le \frac{C}{\mu(B(0,1))} \int_H\exp\Big[
\frac{C\|x\|_H^{2}}{t^{\frac{\sigma+2}{\sigma}}}\Big] \mu(d
x)<\infty$$
holds for sufficiently large $t>0$, since
$\mu(e^{\varepsilon_0\|\cdot\|_H^2})$ is finite according to
 Theorem \ref{T1}$(i)$.

Hence  $P_t$ is hyperbounded operator for sufficient large $t>0$.
 Since $P_t$  has a density w.r.t. $\mu$,  $P_t$
is also compact in $L^2(\mu)$ for large $t>0$ according to \cite[Theorem 2.3]{Wu}. \\


$(iv)$ If $\alpha>2$, then by (\ref{exponential estimate})  we have for small
enough $\varepsilon_0>0$
\beq
d e^{\varepsilon_0\x_H^\alpha}\leq(c-\theta\x_H^{2\alpha-2}
e^{\varepsilon_0\x_H^\alpha})dt+\alpha\varepsilon_0\|X_t\|_H^{\alpha-2}
e^{\varepsilon_0\|X_t\|_H^\alpha}\<X_t, BdW_t \>_H,
\end{equation}
where $c,\theta>0$ are some constants.
By  Jensen's inequality we have
$$ \mathbf{E} e^{\varepsilon_0\x_H^\alpha}\leq
e^{\varepsilon_0\|x\|_H^\alpha}+
ct-\theta\varepsilon_0^{-(2\alpha-2)/\alpha} \int_0^t
\mathbf{E}e^{\varepsilon_0\|X_u\|_H^\alpha} \left(\log
\mathbf{E}e^{\varepsilon_0\|X_u\|_H^\alpha}\right)^{\frac{2\alpha-2}{\alpha}}
du .$$
Let $h(t)$ solve the equation \beq\label{eq1}  h'(t) = c
-\theta\varepsilon_0^{-(2\alpha-2)/\alpha} h(t)\big\{\log
h(t)\big\}^{(2\alpha-2)/\alpha},\ \  \ h(0)=\exp\left[ \varepsilon_0
\left( \|x\|_H^{\alpha}+c\right)\right].
\end{equation}
Then by a standard comparison argument we know \beq\label{eq2}
 \mathbf{E} e^{\varepsilon_0 \|X_t(x)\|_H^{\alpha}}\le
h(t)\le \exp\Big[ c_0\big(1+
t^{-\alpha/(\alpha-2)}\big)\Big],\ \ \ t>0, x\in H \end{equation}
hold for a constant $c_0>0$.
By using (\ref{semigroup inequality}) we have \beq\beg{split}
 \|P_t f\|_\infty &= \|P_{t/2}P_{t/2}
f\|_\infty\\
&\le c_1 \sup_{x\in H} \mathbf{E} \exp \Big[ \frac{c_1}
{t^{(\sigma+2)/\sigma}}\left(1+ \|X_{\frac{t}{2}
}(x)\|_H\right)^{2+\frac{2(2-\alpha)}{\sigma}}\Big],\ \ \ t>0,
\end{split}\end{equation}
where
$c_1>0$ is a constant.
 By  Young's inequality there exists $c_2>0$ such that
$$\frac{c_1} {t^{\frac{\sigma+2}{\sigma}}}
(1+u)^{2+\frac{2(2-\alpha)}{\sigma}}\le \varepsilon_0(1+ u^{\alpha})+ c_2
 t^{-\alpha/(\alpha-2
)},\ \ \ u, t>0.$$
Therefore, there exists a constant   $C>0$ such that
$$\|P_t\|_{L^2(\m)\rightarrow L^\infty(\m)}
 \leq\exp[C(1+t^{-\frac{\aa}{\aa-2}})],\
t>0.$$
The compactness of $P_t$ also follows
 from the \cite{Wu}. \qed

\subsection{Proof of Theorem \ref{T1.3}}

 The proof is based on \cite[Theorem 2.5; 2.6; 2.7]{GM04}.
According to Theorem \ref{T1.2}, we know $\{P_t\}$ is strong Feller and irreducible. Now we only need to
verify the following properties:

(1) For each $r>0$ there exist $t_0>0$ and a  compact set $M\subset H$ such that
$$\inf_{x\in B_r} P_{t_0}\mathbf{1}_M(x)>0, $$
where $B_r=\{y\in H: \ \|y\|_H\le r  \}$.

(2) If $\alpha>2$, then there exist constants $K<\infty$ and $t_1>0$ such that
$$ \mathbf{E}\|X_t(x)\|_H^2\le K,\ x\in H,\ t\ge t_1. $$

(3) If $\alpha=2$, then there exist constants $K<\infty$ and $\beta>0$ such that
$$ \mathbf{E} V\left(X_t(x)\right)\le Ke^{-\beta t}V(x)+K,\ x\in H,\ t\ge 0, $$
where  $V(x)=1+\|x\|_H^2$ and $V(x)=e^{\varepsilon_0\|x\|_H^2}$ for some small constant $\varepsilon_0>0$.

By using the It\^{o} formula we have
$$\|X_t\|_H^2\le \|x\|_H^2+\int_0^t\left(c-\frac{\delta}{2}\|X_s\|_V^\alpha+\gamma\|X_s\|_H^2\right) \d s
+\int_0^t\<X_s, B \d W_s\>_H. $$
If $\alpha>2$, then there exists a constant $c_1>0$
$$\|X_t\|_H^2\le \|x\|_H^2+\int_0^t\left(c_1-\frac{\delta}{4}\|X_s\|_V^\alpha \right) \d s
+\int_0^t\<X_s, B\d W_s\>_H. $$
This implies that there exists $C>0$ such that
\begin{equation}\label{estimate on V-norm}
\mathbf{E}\int_0^t\|X_s\|_V^\alpha \d s \le C(t+\|x\|_H^2),\ t\ge 0.
\end{equation}
And by using  Jensen's inequality
$$\mathbf{E}\|X_t\|_H^2\le \|x\|_H^2+\int_0^t\left[C_1-C_2\left(\mathbf{E}\|X_s\|_H^2\right)^{\alpha/2} \right] \d s. $$
Then by a standard comparison argument we get
$$ \mathbf{E}\|X_t(x)\|_H^2\le C(1+t^{-\frac{2}{\alpha-2}}),\ x\in H,\ t>0.$$
Hence property (2) holds.

According to $(\ref{M})$, for the property (1) it is enough to show that there exist
$t_0$ and  a compact set $M$ in $H$ such that $P_{t_0}\mathbf{1}_M(x)>0$ for some $x\in B_r$.

By (\ref{estimate on V-norm}) and a simple contradiction argument, one can show that  there exists $t_0>0$  such that
$ P_{t_0}\mathbf{1}_M(x)>0 $ for the compact set $M:=\left\{y\in H: \|y\|_V\le \left[C(1+r^2)\right]^{1/\alpha}\right\}$ and $x\in B_r$. So property (1) also holds.

Then the assertions in $(ii)$ hold according to \cite[Theorem 2.5; 2.7]{GM04}. The modified constant in the estimates of spectral gap and exponential convergence comes from the arguments in \cite[Theorem 7.2]{GM06}(in fact, (7.10) implies
that (7.4) holds with that modified constant in \cite{GM06}).

Similarly, if $\alpha=2$ and $\gamma<c_0\delta$, then we can prove
$$  \E\|X_t(x)\|_H^2\le  e^{-\beta t}\|x\|_H^2+ C, \ t\ge 0, x\in H   $$
holds for some constants $\beta>0$ and $C$.  Moreover, by (\ref{exponential estimate})
there  also exists a small constant $\varepsilon_0>0$  such that
$$  \E \exp\left[\varepsilon_0\|X_t(x)\|_H^2\right]\le  e^{-\beta t}e^{\varepsilon_0\|x\|_H^2}+ C, \ t\ge 0, x\in H.   $$
Then the conclusions in $(i)$ follow from  \cite[Theorem 2.5; 2.6]{GM04}.

\section{Application to examples}

To apply our main results, one has to verify condition (\ref{c1})
and (\ref{c2}). To this end, we present some simple sufficient
conditions for (\ref{c1}) and (\ref{c2}).
We first establish the following inequality, which is crucial for verifying (\ref{c1})
in concrete examples.

\beg{lem}\label{L3.1} Let $(E, \langle\cdot,\cdot\rangle)$ be a
 Hilbert space and $\|\cdot\|$ denote its norm, then
for any $r\geq0$ we have \beq\label{3.1} \langle\|a\|^ra-\|b\|^rb,
a-b\rangle\geq2^{-r}\|a-b\|^{r+2}, \ a,b\in E. \end{equation}
\end{lem}
\begin{proof}
We may assume $\|a\|\ge\|b\|$ without loss of
generality. Then we have \beq\beg{split}
& \langle\|a\|^ra-\|b\|^rb,
a-b\rangle\\
=& \|b\|^r\|a-b\|^2+(\|a\|^r-\|b\|^r)\langle a, a-b\rangle\\
=& \|b\|^r\|a-b\|^2+\frac{1}{2}(\|a\|^r-\|b\|^r)(\|a\|^2+\|a-b\|^2-\|b\|^2)\\
\ge& \|b\|^r\|a-b\|^2+\frac{1}{2}(\|a\|^r-\|b\|^r)\|a-b\|^2\\
=& \frac{1}{2}(\|a\|^r+\|b\|^r)\|a-b\|^2\\
\geq& 2^{-r}\|a-b\|^{r+2}. 
\end{split}\end{equation}
\end{proof}

\begin{rem}
  If $r<0$, then (\ref{3.1}) does not hold in general. It's easy to show the
 assumption (\ref{c1}) in Theorem \ref{T1.1} does
  not hold for stochastic fast diffusion equations. For more
  details we refer to \cite{LW}.
\end{rem}

The first example for the application of our main results is the stochastic porous media equation (cf.\cite{PR,L08,W07}).
The main results obtained in \cite{W07} are covered by our main theorems. Moreover, we derive some new results
such as the irreducibility (hence the uniqueness of invariant measures), exponential ergodicity and the existence of spectral gap for stochastic porous media equations here. Now we apply the main results to other types of stochastic evolution equations in Hilbert space.
In the following examples $L(Y,Z)$ denotes the space of all bounded
linear operators from
 $Y$ to $Z$ and $\mathbf{Ran}(B)$ denotes the range of operator $B$.
\beg{exa} (Stochastic reaction-diffusion equation)\\
Let $\Lambda$ be an open bounded domain in $\mathbb{R}^d$ with smooth boundary and
 $\Delta$ be the Laplace operator on $L^2(\Lambda)$ with Dirichlet
boundary condition. Consider
the following triple\\
$$W_0^{1,2}(\Lambda)\cap L^p(\Lambda) \subseteq L^2(\Lambda)\subseteq
 \left(W_0^{1,2}(\Lambda) \cap L^p(\Lambda)  \right)^*$$
and the stochastic reaction-diffusion equation
 \beq\label{srf}
 \d X_t=(\Delta X_t-c|X_t|^{p-2}X_t)\d t+B\d W_t ,\ X_0=x\in L^2(\Lambda)
  \end{equation}
where $ p> 1$ and $c\geq0$,  $B$ is a Hilbert-Schmidt operator
 and $W_t$ is a
 cylindrical Wiener process on $L^2(\Lambda)$, then according to $\cite{Zh}$
$(\ref{srf})$ has a unique strong solution and $(\ref{c1})$ holds with $N(u)=\|u\|_{1,2}^2$ (Sobolev norm).
Hence the assertions in $\mathbf{Theorem\ \ref{T1}}$ hold
for $(\ref{srf})$.

 Moreover,
if $B$  is a one-to-one operator such that
  $$W_0^{1,2}(\Lambda)\subseteq
\mathbf{Ran}(B),\ \
  B^{-1}\in L(W_0^{1,2}(\Lambda);
 L^2(\Lambda)),$$
then $(\ref{c2})$ also holds.
 In particular, if $d=1$ and $B:=(-\Delta)^{-\theta}$
with $\theta\in(\frac{1}{4}, \frac{1}{2}]$,
 then $B$ is a  Hilbert-Schmidt operator and $(\ref{c2})$ holds.
Hence
  the assertions in $\mathbf{Theorem\ \ref{T1.1}, \ref{T1.2}}$ and $\mathbf{\ref{T1.3}}$
 also holds for $(\ref{srf})$.  Particularly, the associated transition semigroup of
$(\ref{srf})$ is hyperbounded and $V$-uniformly ergodic.
\end{exa}

\beg{rem}  Suppose that
$$0<\lambda_1\leq\lambda_2\leq\cdots\leq\lambda_n\leq\cdots$$
are the eigenvalues of $-\Delta$ and the corresponding eigenvectors
$\{e_i\}_{i\geq1}$ form an orthonormal basis on $L^2(\Lambda)$. Assume
$Be_i:=b_ie_i$ and  there exists a positive constant C such that
$$ \sum_i b_i^2<+\infty;\ \   b_i\geq\frac{C}{\sqrt{\lambda_i}}, ~~~ i\geq1, $$
then $B$ is a Hilbert-Schmidt operator and $(\ref{c2})$ holds.

By the Sobolev inequality (\cite{W00},Corollary 1.1 and
3.1) we have
$$\lambda_i \geq c i^{2/d}, \  \  i\geq1$$
hold for some constant $c>0$. This implies that the space
dimension $d$ is less than $2$.
However, if we  consider a general negative
definite self-adjoint operator $L$ instead of $\Delta$ in
(\ref{srf}), e.g. $L:=-(-\Delta)^q, q>0$. Then, by  the
spectral representation theorem, our results can apply to  examples on
 $\mathbb{R}^d$ with $d\ge 2$.  We may refer to \cite{LW,W07} for more details.
\end{rem}

\beg{exa} (Stochastic  $p$-Laplace equation)\\
 Let $\Lambda$ be an open bounded domain in $\mathbb{R}^d$ with smooth boundary. Consider
the  triple
$$W^{1,p}_0(\Lambda)\subseteq L^2(\Lambda)\subseteq
 (W^{1,p}_0(\Lambda))^*$$
and the stochastic $p$-Laplace equation
 \beq\label{sp}
 \d X_t=\left[ \mathbf{div}(|\nabla X_t|^{p-2}\nabla X_t)
 -c|X_t|^{\tilde{p}-2}X_t\right]\d t+B\d W_t , X_0=x,
  \end{equation}
 where $c\geq0$, $2\leq p<\infty, 1\leq\tilde{p}\leq p$,  $B$ is a Hilbert-Schmidt operator
and $W_t$ is
a cylindrical  Wiener process on $L^2(\Lambda)$, then
the assertions in $\mathbf{Theorem\ \ref{T1}}$ hold for $(\ref{sp})$.

Moreover, if $d=1$ and $B:=(-\Delta)^{-\theta}$
with $\theta\in(\frac{1}{4}, \frac{1}{2}]$,
 then $(\ref{c2})$ also holds.
 Therefore the assertions in $\mathbf{Theorem\ \ref{T1.1}, \ref{T1.2}}$ and
 $\mathbf{\ref{T1.3}}$ also hold for $(\ref{sp})$.
 In particular, if $p>2$, then
 the associated transition semigroup of $(\ref{sp})$ is
 ultrabounded and uniformly exponential ergodic, and its generator also has a spectral gap.
\end{exa}

 \begin{proof}
 According to \cite[Example 4.1.9]{PR}, $(A1)-(A4)$ hold for
(\ref{sp}). Hence we only need to
 verify (\ref{c1}) for $N(u)=\|u\|_{1,p}^p$
   under our
 assumptions. By using Lemma \ref{L3.1}
and the Poincar\'{e} inequality we have
 $$\begin{aligned}
 & _{V^*}\langle \mathbf{div}(|\nabla u|^{p-2}\nabla
 u)-\mathbf{div}(|\nabla v|^{p-2}\nabla v), u-v
 \rangle_{V}\\
 &=-\int_\Lambda  \< \ |\nabla u(x)|^{p-2}\nabla u(x)-|\nabla
 v(x)|^{p-2}\nabla v(x), \nabla u(x)-\nabla
 v(x)
 \>_{\mathbb{R}^d} dx\\
 &\leq-2^{p-2}\int_\Lambda |\nabla u(x)-\nabla v(x)|^{p}dx\\
 &\leq -C\|u-v\|_{1,p}^p,\ u,v\in W_0^{1,p}(\Lambda),
\end{aligned}$$ where $C>0$ is a  constant.
And it is  easy to show that
$$ \ _{V^*}\langle
|u|^{\tilde{p}-2}u-|v|^{\tilde{p}-2}v, u-v
 \rangle_{V}\ge 0.$$
Hence (\ref{c1}) holds.

If $d=1$ and $B:=(-\Delta)^{-\theta}$
with $\theta\in(\frac{1}{4}, \frac{1}{2}]$, then
there exists a constant $c>0$ such that (see the remark above)
$$ \|u\|_{1,2} \ge c \|u\|_{B},\ u\in W_0^{1,p}(\Lambda). $$
This implies (\ref{c2}) holds.
\end{proof}

\begin{rem} (1) The Harnack inequality and some consequent properties still hold if
 one  also adds some locally bounded linear (or  order less than $p$) perturbation in the drift.
Only for certain properties (e.g. hyperboundedness or ultraboundedness) we need to require the drift is
dissipative (i.e.$\gamma\le 0$). 

 (2) If we assume $B=0$ in (\ref{sp}), then by  Theorem $\ref{T1}$$(iii)$
we can get the following decay of
the solution to the classical $p$-Laplace equation
$$\sup_{x\in L^2(\Lambda)} \|X_t(x)\|_{L^2}\le C t^{-\frac{1}{p-2}}, \
t>0. $$

\end{rem}

The following SPDE has been studied in \cite{KR,L08}, in which  the main part of
drift in the equation is a high order generalization of the Laplace operator.

\beg{exa}
 Let $\Lambda$ be an open bounded domain in $\mathbb{R}^1$ and
 $m\in\mathbb{N_+}$. Consider
the following triple\\
$$W^{m,p}_0(\Lambda)\subseteq L^2(\Lambda)\subseteq
 (W^{m,p}_0(\Lambda))^*$$
and the stochastic evolution equation
 \beq\label{sh}
 \d X_t(x)=\left[(-1)^{m+1}\frac{\partial^m}{\partial
 x^m}\left(\left|\frac{\partial^m}{\partial x^m} X_t(x)\right|^{p-2}
\frac{\partial^m}{\partial x^m} X_t(x)\right)
 -c|X_t(x)|^{\tilde{p}-2}X_t(x)\right]\d t+B\d W_t ,
  \end{equation}
where $c\geq0$, $2\leq p<\infty, 1\leq\tilde{p}\leq p$, $B\in L_2(L^2(\Lambda))$
 and $W_t$ is
 a cylindrical Wiener process on $L^2(\Lambda)$, then
the assertions in $\mathbf{Theorem \ref{T1}}$ hold for $(\ref{sh})$.

 Moreover,
 if $B$ is also a one-to-one operator such that $B^{-1}\in
 L(W^{m,p}_0(\Lambda);
 L^2(\Lambda))$, then $(\ref{c2})$ is also satisfied.
 Hence the assertions in $\mathbf{Theorem\
 \ref{T1.1}, \ref{T1.2}}$ and $\mathbf{\ref{T1.3}}$ hold for $(\ref{sh})$.
 In particular, the associate transition semigroup of
  the solution is ultrabounded if $p>2$ and hyperbounded if $p=2$.
\end{exa}
\beg{proof}  The proof is similar to the argument in Example 3.3 by
taking $N(u)=\|u\|_{m,p}^p$.
\end{proof}

\begin{rem}(i)
 If we assume $p>2$ and $B=0$ in (\ref{sh}), then by Theorem \ref{T1}
we also  obtain the decay of
the solution to the  deterministic evolution equation, i.e.
$$\sup_{f\in L^2(\Lambda)} \|X_t^f\|_{L^2}\le C t^{-\frac{1}{p-2}}, \
t>0, $$
where $X_t^f$ denote the solution to the following equation
$$\frac{\d X_t(x)}{\d t}=(-1)^{m+1}\frac{\partial^m}{\partial
 x^m}\left(\left|\frac{\partial^m}{\partial x^m} X_t(x)\right|^{p-2}
\frac{\partial^m}{\partial x^m} X_t(x)\right)
 -c|X_t(x)|^{\tilde{p}-2}X_t(x), \ X_0=f\in L^2(\Lambda).$$

(ii) Assume that
$$0<\lambda_1\leq\lambda_2\leq\cdots\leq\lambda_n\leq\cdots$$
are the eigenvalues of a positive definite self-adjoint operator $L$
where $\mathcal{D}(\sqrt{L})=W^{m,2}_0(\Lambda)$, the corresponding
eigenvector $\{e_i\}_{i\geq1}$ is an ONB of $L^2(\Lambda)$. Suppose
$Be_i:=b_ie_i$ and  there exists a  constant $C>0$ such
that
$$\sum_i b_i^2<+\infty;\ \   b_i\geq\frac{C}{\sqrt{\lambda_i}},  ~~~ i\geq1, $$
then $B$ is a Hilbert-Schmidt operator on $L^2(\Lambda)$ and (\ref{c2}) is satisfied.

\end{rem}

\section*{Acknowledgements}
The author would like to thank Professor Michael R\"{o}ckner, Professor Fengyu Wang
 and Professor Bohdan Maslowski
for their valuable discussions and suggestions. The useful comments
from the referee is also gratefully  acknowledged.

\end{document}